\documentclass[11pt]{article}
\usepackage{latexsym,amssymb}
\usepackage{amsmath}
\usepackage{color}
\usepackage{mathrsfs}
\usepackage[abbrev]{amsrefs}

\newtheorem{thm}{Theorem}[section]

\newtheorem{prop}[thm]{Proposition}
\newtheorem{lem}[thm]{Lemma}

\numberwithin{equation}{section}
\numberwithin{thm}{section}

\begin{document}

\title{ 
The semigroup generated by the Dirichlet Laplacian 
\\
of fractional order}
\author{
Tsukasa Iwabuchi\\
Mathematical Institute\\
Tohoku University\\
Sendai 980-8578 Japan
}
\date{}
\maketitle
\footnote[0]
{


E-mail: t-iwabuchi@m.tohoku.ac.jp

}

\begin{abstract}
In the whole space $\mathbb R^d$, 
linear estimates for heat semi-group in 
Besov spaces are well established, 
which are estimates of $L^p$-$L^q$ type, maximal regularity, 
e.t.c. 
This paper is concerned with such estimates for semi-group generated by 
the Dirichlet Laplacian of fractional order 
in terms of the Besov spaces on an arbitrary open set of $\mathbb R^d$. 
\end{abstract}
\section{Introduction}

\quad 
Let $\Omega$ be an arbitrary open set of $\mathbb R^d$ with $ d \geq 1$. 
We consider the Dirichlet Laplacian $A$ on $L^2 (\Omega)$, namely, 
\begin{equation}\notag 
A  = -\Delta 
= - \sum _{j = 1} ^d \frac{\partial ^2}{\partial x_j ^2} 
\end{equation}
with the domain 
$$
\mathcal D (A) 
:= 
\{ \, f \in H^1_0 (\Omega ) 
  \, | \, 
  \Delta f \in L^2 (\Omega) 
  \, \}. 
$$
We consider the fractional Laplacian and semi-group 
\begin{equation}\notag 
A^{\frac{\alpha}{2}}
= \int_{-\infty}^\infty \lambda ^{\frac{\alpha}{2}} dE_A (\lambda), 
\qquad 
e^{-t A ^{\frac{\alpha}{2}}} 
= \int_{-\infty} ^\infty e^{-t \lambda ^{\frac{\alpha}{2}}} d E_A (\lambda) , 
\quad t \geq 0 . 
\end{equation}
Here, $\alpha > 0$ and $\{ E_A (\lambda) \}_{ \lambda \in \mathbb R}$ denotes 
the spectral resolution of identity, where it is determined uniquely for 
the self-adjoint operator $A$ by the spectral theorem. 
The motivation of study of fractional Laplacian comes from 
the study of fluid mechanics, stochastic process, finace, etc., 
see for instance \cite{App_2009,Be_1996,CKS-2010,VIKH-2009} and references therein. 
We also refer to \cite{DiPaVa-2012,Vaz-2012,Vaz-2014} 
where one can find some results on fundamental properties of fractional Sobolev spaces 
and applications to partial differential equations. 
\\

In the paper \cite{IMT-preprint2}, 
based on the spectral theory for the Dirichlet Laplacian $A$ on $L^2 (\Omega) $, 
a kind of $L^p$ theory was established and the Besov spaces on an open set $\Omega$ were introduced, 
where regularity of functions is measured by $A$. 
The purpose of this paper is to develop the linear estimates 
for the semi-group generated by the Dirichlet Laplacian of fractional order 
in the homogeneous Besov spaces $\dot B^s_{p,q} (A)$, namely, 
the estimate of $L^p$-$L^q$ type, smoothing effects, 
continuity in time of semi-group, 
equivalent norms with the semi-group 
and maximal regularity estimates. 
Such estimates with heat semi-group in the case when $\Omega = \mathbb R^d$ 
are well established
(see \cite{BaChDa_2011,Che-2004,D-2005,D-2007,DM-2009,HP-1997,KOT-2003,Lem_2002,OS-2009,OS-pre,Pee_1976}). 
In this paper we consider open sets of $\mathbb R^d$ and 
the semi-group generated by the fractional Laplacian 
with the Dirichlet boundary condition. 
\\

Let us recall the definitions of spaces of test functions and tempered distributions and the Besov spaces 
associated with 
the Dirichlet Laplacian (see~\cite{IMT-preprint2}). 
We take $\phi_0(\cdot) \in C^\infty_0(\mathbb R)$ a non-negative function on $\mathbb R$ such that  
\begin{equation}
\label{325-1}
{\rm supp \, } \phi _0
\subset \{ \, \lambda \in \mathbb R \, | \, 2^{-1} \leq \lambda \leq 2 \, \}, 
\quad \sum _{ j \in \mathbb Z} \phi_0 ( 2^{-j}\lambda) 
 = 1 
 \quad \text{for } \lambda > 0, 
\end{equation}
and $ \{ \phi_j \}_{j \in \mathbb Z}$ is defined by letting 
\begin{gather} \label{325-2}
\phi_j (\lambda) := \phi_0 (2^{-j} \lambda) 
 \quad \text{for }  \lambda \in \mathbb R . 
\end{gather}

\vskip3mm 

\noindent 
{\bf Definition. } 
{\it 
{\rm (i)} 
{\rm (}Linear topological spaces
$\mathcal X_0 (\Omega)$ and $\mathcal X^\prime_0 (\Omega)${\rm ).}
$\mathcal X_0 (\Omega)$ is 
defined by letting
\begin{equation}\notag 
\mathcal X_0 (\Omega) 
:= \big\{ f \in  L^1 (\Omega) \cap \mathcal D (A) 
 \, \big| \, 
    A^{M} f \in L^1(\Omega ) \cap \mathcal D (A) \text{ for all } M \in \mathbb N 
   \big\} 
\end{equation} 
equipped with the family of semi-norms $\{ p_{0,M}
 (\cdot) \}_{ M = 1 } ^\infty$ 
given by 
\begin{equation}\notag 
p_{0,M}(f) := 
\| f \|_{ L^1(\Omega)} 
+ \sup _{j \in \mathbb N} 2^{Mj} 
  \| \phi_j (\sqrt{A}) f \|_{ L^1(\Omega)} . 
\end{equation}

\noindent 
{\rm (ii)}
{\rm (}Linear topological spaces
$\mathcal Z_0 (\Omega)$ and $\mathcal Z^\prime_0 (\Omega)${\rm ).} 
$\mathcal Z_0 (\Omega)$ is 
defined by letting
\begin{equation}\notag 
\mathcal Z_0 (\Omega) 
:= \Big\{ f \in \mathcal X_0 (\Omega) 
 \, \Big| \, 
  \sup_{j \leq 0} 2^{ M |j|} 
    \big\| \phi_j \big(\sqrt{ A } \big ) f \big \|_{L^1(\Omega)} < \infty 
  \text{ for all } M \in \mathbb N
   \Big\} 
\end{equation}
equipped 
with the family of semi-norms $\{ q_{0,M} (\cdot) \}_{ M = 1}^\infty$ 
given by 
\begin{equation}\notag 
q_{0,M}(f) := 
\| f \|_{L^1 (\Omega) }
+ \sup_{j \in \mathbb Z} 2^{M|j|} \| \phi_j (\sqrt{A}) f \|_{L^1(\Omega)}. 
\end{equation}
}   

\noindent 
{\bf Definition. } 
{\it 
For $s \in \mathbb R$ and $1 \leq p,q \leq \infty$, 
$\dot B^s_{p,q} (A) $ is defined by letting
\begin{equation} \notag 
\dot B^s_{p,q}  (A)
:= \{ f \in \mathcal Z'_0(\Omega) 
     \, | \, 
     \| f \|_{\dot B^s_{p,q}(A)} < \infty 
   \} , 
\end{equation}
where
\begin{equation}\notag 
     \| f \|_{\dot B^s_{p,q}(A)} 
       := \big\| \big\{ 2^{sj} \| \phi_j (\sqrt{A}) f\|_{L^p(\Omega)}
                 \big\}_{j \in \mathbb Z}
          \big\|_{\ell^q (\mathbb Z)}. 
\end{equation}
}   

Let us mention the basic 
properties of $\mathcal X_0 (\Omega), 
\mathcal Z_0 (\Omega)$, their duals, $\dot B^s_{p,q} (A)$ and explain the operators 
$\phi_j (\sqrt{A})$ and the Laplacian of fractional order. 

\vskip3mm

\noindent 
{\bf Proposition}. (\cite{IMT-preprint2})
{\it Let $s ,\alpha\in \mathbb R$, $1\leq p,q,r \leq \infty$. 
Then the following hold: 
\begin{itemize}
\item[{\rm (i)}] 
 $\mathcal X_0 (\Omega)$ and $ \mathcal Z_0 (\Omega)$ 
are Fr\'echet spaces and enjoy 
$
\mathcal X_0 (\Omega) \hookrightarrow L^p (\Omega) 
\hookrightarrow \mathcal X_0 ' (\Omega), 
\mathcal Z_0 (\Omega) \hookrightarrow L^ p (\Omega)
\hookrightarrow \mathcal Z_0 ' (\Omega). 
$

\item[{\rm (ii)}] 
$\dot B^s_{p,q} (A)$ is a Banach space and enjoys 
$
\mathcal Z_0 (\Omega) \hookrightarrow \dot B^s_{p,q} (A) 
\hookrightarrow \mathcal Z_0 ' (\Omega). 
$

\item[{\rm (iii)}]
If $p,q < \infty$ and 
$1/p + 1/p' = 1/q + 1/q' = 1$, 
the dual space of $\dot B^s_{p,q} (A)$ is 
$\dot B^{-s}_{p',q'} (A)$. 

\item[{\rm (iv)}] 
If $r \leq p$, 
$\dot B^{s+d(\frac{1}{r}-\frac{1}{p})}_{r,q} (A)$ 
is embedded to $\dot B^s_{p,q} (A)$. 

\item[{\rm (v)}] 
For any $f \in \dot B^{s+\alpha}_{p,q} (A)$,  
$A^{\frac{\alpha}{2}} f \in \dot B^{s}_{p,q} (A) $. 
\end{itemize}
} 

\vskip3mm

It should be noted that $\phi_j (\sqrt{A})$ and $A$ 
are defined in $L^2 (\Omega)$ initially and  
by the argument in \cite{IMT-preprint2} they can be realized 
as operators in $\mathcal Z_0 '(\Omega)$ and Besov spaces. 
In the proof, the uniform boundedness in $L^p (\Omega) $ 
of $\phi_j (\sqrt{A})$ with respect to $j \in \mathbb Z$ is essential. 
Uniformity in $L^2 (\Omega)$ is proved easily by the spectral theorem, 
while that in $L^1 (\Omega)$ is not trivial. 
For any open set $\Omega \subset \mathbb R ^d$, 
$L^1 (\Omega)$ boundedness is known 
in some papers (see Proposition~6.1 in \cite{DOS-2002} 
and also Theorem~1.1 in \cite{IMT-RMI}). 
Let us explain the strategy of its proof 
along \cite{IMT-RMI}  
(see also a comment below Lemma~\ref{lem:317-2}). 
The uniform boundedness in $L^1 (\Omega)$ is proved via estimates in 
amalgam spaces $\ell^1 (L^2)_{\theta} $, where the side length of 
each cube is scaled by $\theta ^{\frac{1}{2}}$, $\theta = 2^{-2j}$ (see section 2), 
together with the Gaussian upper bounds 
of the kernel of $e^{-tA}$. 
That scaling fits for the scaled 
operator $\phi_j (\sqrt{A}) = \phi _0 (2^{-j} \sqrt{A})$, 
and we can handle the norm in $\ell ^1 (L^2) _ \theta$ 
through the estimates in $L^2 (\Omega)$, since its norm is defined 
locally with $L^2 (\Omega)$. The Gausian upper bounds 
of the kernel of $e^{-tA}$ is necessary in order 
to estimate $L^1 (\Omega)$ norm via $\ell ^1 (L^2)_\theta $. 
Once $L^1 (\Omega)$ estimate is proved, $L^p (\Omega)$ case 
is assured by the interpolation and the duality argument. 

\vskip3mm

As to the Laplacian of fractional order, it was shown 
in the proof of Proposition~3.2 in \cite{IMT-preprint2} that 
$A^{\frac{\alpha}{2}}$ is a continuous operator from 
$\mathcal Z_0 ' (\Omega)$ to itself, 
which is proved as follows: Show the continuity of $A^{\frac{\alpha}{2}}$ in 
$\mathcal Z_0 (\Omega)$ first with 
the boundedness of spectral multipliers 
$$
\| A^{\frac{\alpha}{2}}\phi_j (\sqrt{A}) 
\|_{L^1 (\Omega) \to L^1 (\Omega) } 
\leq C 2^{\alpha j}
$$
for all $j \in \mathbb Z$ and consider their dual operator 
together with the approximation of the identity
$$
f = \sum _{ j \in \mathbb Z} \phi_j (\sqrt{A}) f
\quad \text{in } \mathcal Z_0 ' (\Omega) 
\quad \text{for any } f \in \mathcal Z_0 ' (\Omega) . 
$$
Hence, we define $A^{\frac{\alpha}{2}}$ by 
\begin{equation}\notag
A^{\frac{\alpha}{2}} f = 
\sum _{ j \in \mathbb Z} \big( A^{\frac{\alpha}{2}}\phi_j (\sqrt{A}) \big) f 
\quad \text{in } \mathcal Z_0 ' (\Omega)
\quad \text{for any  } f \in \mathcal Z_0 ' (\Omega) .
\end{equation}
Noting that $e^{-t A^{\frac{\alpha}{2}}} \phi_j (\sqrt{A})$ with $t \geq 0$ 
is also bounded in $L^1 (\Omega)$ (see Lemma~\ref{lem:317-1} and 
\eqref{317-13} below), 
we also define $e^{-t A^{\frac{\alpha}{2}}}$ by 
$$
e^{-t A^{\frac{\alpha}{2}}} f = 
\sum _{ j \in \mathbb Z} \big( e^{-t A^{\frac{\alpha}{2}}} \phi_j (\sqrt{A}) \big) f 
\quad \text{in } \mathcal Z_0 ' (\Omega)
\quad \text{for any  } f \in \mathcal Z _0 ' (\Omega) .
$$

\vskip3mm

We state four theorems on the semi-group generated by $A^{\frac{\alpha}{2}}$;  
the estimates of $L^p$-$L^q$ type and smoothing effects, 
continuity in time, equivalent norms with semi-group  
and maximal regularity estimates, 
referring to the results in the case when 
$\Omega = \mathbb R^d $ and $\alpha = 2$. 
\\

We start by considering estimates of $L^p$-$L^q$ type and smoothing effects. When $\Omega = \mathbb R^n$, 
it is well known that 
$$
\| e^{t \Delta} f\|_{L^q (\mathbb R^d)} 
\leq C t^{-\frac{d}{2} (\frac{1}{p}-\frac{1}{q}) } \| f \|_{L^p (\mathbb R^d)}, 
\quad 
\| \nabla e^{t\Delta} f \|_{L^p (\mathbb R^d)}  
\leq C t^{-\frac{1}{2}} \| f \|_{L^p (\mathbb R^d)}, 
$$
where $1 \leq p,q \leq \infty$, $f \in L^p (\mathbb R^d)$. 
Hence one can show that 
$$
\big\| e^{ t \Delta } f \big\|_{\dot B^{s_2}_{p_2,q}(A)} 
\leq C t^{-\frac{d}{2} (\frac{1}{p_1} - \frac{1}{p_2}) 
          - \frac{s_2 - s_1}{ 2}
         } 
     \| f \|_{\dot B^{s_1}_{p_1 , q} (A) }, 
$$
where $s_2 \geq s_1$, $1 \leq p_1 \leq p_2 \leq \infty$ and $1 \leq q \leq \infty$. 
The following is the linear estimates for the semi-group 
generated by $A^\frac{\alpha}{2}$ on an open set.

\begin{thm}\label{thm:bdd}
Let $\alpha  > 0$, $ t \geq 0$, 
$s, s_1 ,s_2 \in \mathbb R$, $1 \leq p,p_1,p_2, q, q_1 , q_2 \leq \infty$. 
\\
{\rm (i)} 
$e^{-tA^{ \frac{\alpha }{2} }}$ is a bounded linear operator in 
$\dot B^s_{p,q } (A)$, i.e., 
there exists a constant $C > 0$ such that 
 for any $f \in \dot B^s_{p,q} (A)$
\begin{equation}\label{317-11}
 e^{-tA^{ \frac{\alpha }{2} }} f \in \dot B^s_{p,q} (A)
\quad \text{and} \quad 
\big\| e^{-tA^{ \frac{\alpha }{2} }} f \big\|_{\dot B^{s}_{p,q}(A)} 
\leq C \| f \|_{\dot B^{s}_{p , q} (A) }. 
\end{equation} 
\\
{\rm (ii)} If $s _2 \geq s_1$, $p_1 \leq p_2$ and 
$$
d \Big(\frac{1}{p_1} - \frac{1}{p_2} \Big) 
          + s_2 - s_1 > 0, 
$$
then there exists a constant $C > 0$ such that 
\begin{equation}\label{317-12}
\big\| e^{-tA^{ \frac{\alpha }{2} }} f \big\|_{\dot B^{s_2}_{p_2,q_2}(A)} 
\leq C t^{-\frac{d}{\alpha} (\frac{1}{p_1} - \frac{1}{p_2}) 
          - \frac{s_2 - s_1}{\alpha}
         } 
     \| f \|_{\dot B^{s_1}_{p_1 , q_1} (A) } 
\end{equation}
for any $f \in \dot B^{s_1}_{p_1, q_1} (A)$. 
\end{thm}
{\bf Remark. } 
On the estimate \eqref{317-12}, the regularity on indices $q_1$ and $q_2$ 
is gained without loss of the singularity at $t = 0$. This estimate is 
known in the case when $\Omega = \mathbb R^n $ and $\alpha  = 2$ (see \cite{KOT-2003}).
\\

As to the continuity in time of heat semi-group $e^{t \Delta }$ 
when $\Omega = \mathbb R^d$, 
it is well known that for $1 \leq p < \infty$ 
\begin{equation}\notag 
\lim _{t \to 0} \| e^{t\Delta} f - f \|_{L^p (\mathbb R^d)} = 0 
\quad \text{for any } f \in L^ p (\mathbb R^d).
\end{equation}
In the case when $p = \infty$, the above strong convergence does not hold 
in general while it holds in the dual weak sense. 
The following theorem is concerned with such continuity 
in the Besov spaces on an open set.

\begin{thm}\label{thm:cont}
Let $s \in \mathbb R$, 
$1 \leq p,q \leq \infty$ and 
$1/p + 1/p' = 1/q + 1/q' = 1$. 
\\
{\rm (i)} Assume that $q < \infty$ and $f \in \dot B^s_{p,q} (A)$. 
Then 
\begin{equation}\notag 
\lim _{t \to 0} \big\| e^{-tA^{ \frac{\alpha }{2} }} f -f 
\big\|_{\dot B^{s}_{p,q}(A)} = 0.
\end{equation}
{\rm (ii)} 
Assume that $1 < p \leq \infty$, $q = \infty$ and $f \in \dot B^s_{p,\infty} (A)$. 
Then $e^{-t A^{\frac{\alpha}{2}}} f$ converges to $f$ in the dual weak sense 
as $t \to 0$, namely, 
$$
\lim _{ t \to 0 } 
\sum _{j \in \mathbb Z} 
\int _{\Omega} 
   \Big\{ \phi_j (\sqrt{A})\big( e^{-tA^{ \frac{\alpha }{2} }} f -f \big)\, 
   \Big\}
   \overline{ \Phi_j (\sqrt{A})g} \, dx 
= 0 
$$
for any $g \in \dot B^{-s}_{p',1} (A)$. 

\end{thm}

\noindent 
{\bf Remark. } 
Rerated to Theorem~\ref{thm:cont} {\rm (ii)}, 
it should be noted that the predual of $\dot B^{s}_{p,q} (A)$ 
is $\dot B^{-s}_{p',q'} (A)$ for $1 < p,q \leq \infty$, 
where $1/p + 1/p' = 1/q + 1/q' = 1$. 
In fact, we can regard $f \in \dot B^s_{p,q} (A)$ as an element of 
the dual of $\dot B^{-s}_{p',q'} (A)$ by 
$$
\langle f , g \rangle 
= \sum _{j \in \mathbb Z} 
\int _{\Omega} 
   \big\{ \phi_j (\sqrt{A})f  
   \big\}
   \overline{ \Phi_j (\sqrt{A})g} \, dx 
$$
for any $g \in \dot B^{-s}_{p',q'} (A)$ (see \cite{IMT-preprint2}), 
where $\Phi_j := \phi_{j-1} + \phi_j + \phi_{j+1}$. 
\\

As to the characterization of norm by using 
semi-group when $\Omega = \mathbb R^d$, it is 
known that 
$$
\| f \|_{\dot B^s_{p,q} (A)} \simeq 
\Big\{ \int _0^\infty 
     \Big(  t^{-\frac{s}{2}} \|  e^{t \Delta } f \|_{L^p (\mathbb R^d)} 
     \Big)^q
     \frac{dt}{t}
\Big\} ^{\frac{1}{q}} 
$$
where $s < 0$ (see e.g. \cite{Lem_2002}). 
We consider the equivalent norm of Besov spaces on an open set 
by using the semi-group generated by $A^\frac{\alpha}{2}$.

\begin{thm}\label{thm:equiv}
Let $\alpha > 0$, $s , s_0\in \mathbb R$, 
$s_0 > s/ \alpha$ and $1 \leq p,q \leq \infty$. 
Then there exists a constant $C>0$ such that 
\begin{equation}\label{318-1}
C^{-1}\| f \|_{\dot B^s_{p,q} (A)} 
\leq 
\Big\{ \int _0^\infty 
     \Big(  t^{-\frac{s}{\alpha}} \| (tA^{\frac{\alpha}{2}} ) ^{s_0} e^{-tA^{\frac{\alpha}{2}}} f \|_{X} 
     \Big)^q
     \frac{dt}{t}
\Big\} ^{\frac{1}{q}}
\leq C \| f \|_{\dot B^s_{p,q} (A)} 
\end{equation}
for any $f \in \dot B^s_{p,q} (A)$, 
where $X = L^p (\Omega), \dot B^0_{p,r} (A)$ with $1 \leq r \leq \infty$. 
\end{thm}

Since the equivalence \eqref{318-1} is closely related to the real interpolation in the Besov spaces, 
we also mention that the interpolation is also available 
(see e.g. \cite{BL_1976,Triebel_1983} and also Proposition \ref{prop:real} 
in Appendix A below). 
\\

The last result is concerned with the maximal regularity estimates. 
When $\Omega = \mathbb R^d$, 
the Cauchy problem which we should consider is the following,  
$$
\begin{cases}
\partial_t u - \Delta u = f ,
& \quad t > 0, \, x \in \mathbb R^d,
\\
u(0,x ) = u_0 (x), 
& \quad x \in \mathbb R^d. 
\end{cases}
$$
For $1 < p,q < \infty$, 
the solution $u$ of the above problem satisfies that 
$$
\| \partial _t u \|_{L^q (0,\infty ; L^p (\mathbb R^d))}
+ \| \Delta u \|_{L^q (0,\infty ; L^p (\mathbb R^d))} 
\leq 
C\| u_0 \|_{\dot B^{ 2 - \frac{2}{q}}_{p,q} (A) }
+ C \| f \| _{L^q (0,\infty ; L^p (\mathbb R^d))} ,
$$
provided that $u_0\in \dot B^{ 2 - \frac{2}{q}}_{p,q} (A) $ and 
$f \in L^q (0,\infty ; L^p (\mathbb R^d))$ (see \cite{HP-1997,Lem_2002}). 
We note that maximal regularity 
such as the above is well studied in the general framework on Banach spaces with 
unconditional martingale differences which is called UMD 
(see~\cite{Am_1995,DaGr-1975,DeHiPr_2003,DoVe-1987,Lad_1968,Weis-2001}), and 
that the cases when $p,q = 1,\infty$ require a different treatment 
from UMD since the spaces are not reflexive. 
In terms of Besov spaces, 
one can consider $\dot B^0_{p,q} (A)$ 
for all indices $p,q$ with $1 \leq p,q \leq \infty$ 
(see \cite{D-2005,D-2007,DM-2009,HP-1997,OS-2009,OS-pre}). 
Our result on the maximal regularity estimates on open set 
is formulated in the following way.

\begin{thm}\label{thm:maxi}
Let $s \in \mathbb R$, $\alpha > 0$ and $1 \leq p,q \leq \infty$. 
Assume that $u_0 \in \dot B^{s+\alpha - \frac{\alpha}{q}}_{p,q} (A) $, 
$ f \in L^q (0,\infty ; \dot B^s_{p,q}(A))$. Let $u $ be given by  
\begin{gather}
\notag 
u(t) = e^{- t A ^{\frac{\alpha}{2}} } u_0 
   + \int_0^t e^{-(t-\tau) A^{\frac{\alpha}{2}}} f(\tau) d\tau . 
\end{gather}
Then there exists a constant $C>0$ independent of $u_0$ and $f$ such that 
\begin{equation}\label{318-2}
\| \partial _t u \|_{L^q (0,\infty ; \dot B^s_{p,q}(A))}
+ \| A^{\frac{\alpha}{2}} u \|_{L^q (0,\infty ; \dot B^s_{p,q}(A))} 
\leq 
C\| u_0 \|_{\dot B^{s+\alpha - \frac{\alpha}{q}}_{p,q} (A) }
+ C \| f \| _{L^q (0,\infty ; \dot B^s_{p,q}(A))} . 
\end{equation}
\end{thm}


The proof of theorems is based on the boundedness of 
spectral multiplier of the operator 
$ e^{-t A^{\frac{\alpha}{2}}} \phi_j (\sqrt{A})$: 
$$
\| e^{-t A^{\frac{\alpha}{2}}} \phi_j (\sqrt{A}) 
\|_{L^p (\Omega) \to L^p (\Omega)} 
\leq C 
\| e^{-t 2^{\alpha j}(\cdot)^\alpha} \phi_0 ( \sqrt{\cdot} \, ) \|_{H^s (\mathbb R)} 
\quad \text{for all } j \in \mathbb Z, 
$$
where $s > d/2 + 1/2$ (see Lemma \ref{lem:317-1} below). 
The above inequality implies that 
$$
\| e^{-t A^{\frac{\alpha}{2}}} \phi_j (\sqrt{A}) 
\|_{L^p (\Omega) \to L^p (\Omega)} 
\leq C e^{-C^{-1} t 2^{\alpha j}},
$$
and this estimate allows us to show our theorems 
in the analogous method to those in the case when $\Omega = \mathbb R^d$. 
In this paper, we give proofs of all theorems 
by estimating directly for the paper to be self-contained. 
Here, we note that our proofs can be applicable to the estimates for $e^{-tA}$ 
in the inhomogeneous 
Besov spaces and hence similar theorems are able to be obtained. 
On the other hand, for the semi-group generated by 
the fractional Laplacian, since there appears to be problems around low frequency, 
we show only the result for the heat semi-group in Section 7 
(see Theorem~\ref{thm:inhom} below). 
It should be also noted that our argument can be applied 
not only to the Dirichlet Laplacian but also to more general operators 
$A$ such that the Gaussian upper bounds for $e^{-tA}$ hold. 
\\

This paper is organized as follows. In Section 2, 
we prepare a lemma to prove our results. 
Sections 3--6 are devoted to proofs of theorems. 
In Section 7, we state the result for the inhomogeneous Besov spaces. 
In Appendix A, we show the characterization of 
Besov spaces by real interpolation. 
\\

Before closing this section, we introduce some notations. 
We denote by 
$\| \cdot \|_{L^p}$ the $L^p (\Omega)$ norm and 
$\| \cdot \|_{\dot B^s_{p,q}}$ the $\dot B^s_{p,q} (A)$ norm. 
We use the notation $\| \cdot \|_{H^s(\mathbb R)}$ 
as the $H^s (\mathbb R)  $ norm for functions, e.g. 
$\phi_j (\lambda)$, $e^{-t\lambda ^{\frac{\alpha}{2}}}$, 
whose variables are spectral parameter. 
We denote by $\mathcal S (\mathbb R)$ the Schwartz class.

\section{Preliminary}

\quad In this section we introduce the following lemma 
on the boundedness of the scaled spectral multiplier.

\begin{lem}\label{lem:317-1}
Let $N > d/2$, $1 \leq p \leq \infty$, $\delta > 0$ and $a,b > 0$. 
Then there exists a constant $C>0 $ such that 
for any $\phi \in C_0 ^\infty (\mathbb R)$ with  
${\rm supp \, } \phi \subset [a,b] $, 
$G \in C^\infty((0,\infty )) \cap C (\mathbb R)$
and $f \in L^p (\Omega)$ 
\begin{equation}\label{317-1}
\| G(\sqrt{A}) \phi ( 2^{-j} \sqrt{A})  f \|_{L^p} 
\leq C \| G( 2^j \sqrt{\cdot } \, ) \phi (\sqrt{\cdot} \, ) 
       \|_{H^{N+\frac{1}{2} + \delta} (\mathbb R)} \| f \|_{L^p}
\end{equation}
for all $j \in \mathbb Z$. 
\end{lem}

\noindent 
{\bf Remark. } As is seen from the proof below, 
the constant $C$ in the right member of \eqref{317-1} depends on the interval 
$[a,b]$ containing the support of $\phi$. 
\\

To prove Lemma~\ref{lem:317-1}, 
we introduce a set $\mathscr A_N$ of some bounded operators on $L^2 (\Omega)$ 
and scaled amalgam spaces $\ell ^1 (L^2)_{\theta}$ for $\theta > 0$ 
to prepare a lemma. 
Hereafter, for $k \in \mathbb Z^d$, 
$C_\theta (k)$ denotes a cube with the center $\theta^{\frac{1}{2}} k$ 
and side length $\theta ^{\frac{1}{2}}$, namely, 
$$
C_{\theta} (k) 
:= \big\{ x \in \Omega \, \big| \,
    |x_j - \theta^{\frac{1}{2}} k_j| 
    \leq 2^{-1} \theta^{\frac{1}{2}} \text{ for } j = 1,2, \cdots, d
    \, \big\}, 
$$ 
and $\chi _{C_\theta (k)}$ is a characteristic function 
whose support is $C_\theta (k)$. 
\\

\noindent 
{\bf Definition. }  
For $N \in \mathbb N$, 
$\mathscr{A}_N$ denotes the set of all 
bounded operators $T$ on $L^2 (\Omega)$ such that 
$$
\| T \|_{\mathscr{A}_N} := 
\sup_{k \in \mathbb{Z}^d} 
 \big\| |\cdot - \theta^{1/2}k|^{N} T \chi_{C_{\theta}(k)}
 \big\|_{L^2 \to L^2} < \infty. 
 $$

\vskip3mm 

\noindent 
{\bf Definition. } 
The space $\ell ^1 (L^2)_{\theta}$ is defined by letting  
$$
\ell ^1 (L^2)_{\theta} := 
\big\{ f \in L^2 _{\rm loc} (\overline \Omega) \, \big| \, 
   \| f \|_{\ell^1 (L^2)_{\theta}} < \infty
  \, \big\}, 
$$
where
$$
\displaystyle 
\| f \|_{\ell^1 (L^2)_{\theta}} 
:= \sum _{k \in \mathbb Z^d} 
   \| f \|_{L^2 (C_\theta (k))}. 
$$

\begin{lem}\label{lem:317-2} {\rm (\cite{IMT-RMI,IMT-preprint})}
\noindent 
{\rm (i)} 
Let $N \in \mathbb N$ and $N > d/2$. Then 
there exists a constant $C > 0$ such that 
\begin{equation}
\label{317-6}
\| T \|_{\ell^1 (L^2)_{\theta} \to \ell^1 (L^2)_{\theta}} 
\leq C 
  \Big( \| T \|_{L^2 \to L^2} 
    + \theta ^{-\frac{d}{4}} 
      \| T \|_{\mathscr A _N} ^\frac{d}{2N} 
      \| T \|_{L^2 \to L^2} ^{ 1- \frac{d}{2N}}
  \Big) 
\end{equation}
for any $T \in \mathscr A_N$ and $\theta >0 $. 

\noindent 
{\rm (ii)} 
Let $N \in \mathbb N $.   
Then there exists a constant $C > 0$ such that 
\begin{equation}
\label{317-7}
\| \psi ( (M + \theta A)^{-1} ) \|_{\mathscr A_N} 
\leq C \theta ^{\frac{N}{2}} 
   \int_{-\infty}^\infty (1+ |\xi|^2)^{\frac{N}{2}} | \widehat \psi (\xi) | d\xi
\end{equation}
for any $\psi \in \mathcal S (\mathbb R)$ and $\theta > 0$. 
\\
\noindent 
{\rm (iii)} 
Let $M > 0 $ and $\beta > d /4$. Then there exists a constant $C > 0$ such that 
\begin{equation}
\label{317-8}
\big\| (M + \theta A) ^{-\beta} \big\| _{L^1 \to \ell^1 (L^2)_{\theta}} 
\leq C \theta^{ -\frac{d}{2} } 
\end{equation}
for any $\theta > 0$.

\end{lem}

\vskip3mm

\noindent 
{\bf Remark}. 
Lemma~\ref{lem:317-2} is useful to prove $L^1$ boundedness 
of spectral multipliers 
and let us remind how to prove Lemma~\ref{lem:317-2} along 
\cite{IMT-RMI,IMT-preprint} briefly. 
The original idea is from the papers by Jensen-Nakamura~\cite{JN-1994,JN-1995}, 
who studied the Schr\"odinger operators on $\mathbb R^d$. 
On the first inequality \eqref{317-6}, 
we start by decomposing into 
$T = \sum _{m, k \in \mathbb Z^d} \chi_{C_\theta (m)} T \chi _{C_\theta (k)}$, 
and it suffices to show that 
for each $k \in \mathbb Z$ a sum of operator norms 
$
 \sum _{ m \in \mathbb Z} 
 \| \chi_{C_\theta (m)} T \chi_{C_\theta (k)} \|_{L^2 \to L^2}
$
is bounded by the right member of \eqref{317-6}. 
The first term $\| T \|_{L^2 \to L^2}$ 
is obtained just by applying $L^2 (\Omega)$ boundedness 
to $L^2 (C_\theta (m))$ norm with $m = k$. 
The second term is obtained by decomposing the sum into two cases 
when $0 <|m -k| \leq \omega$ and $|m-k| > \omega$ for $\omega > 0$, 
applying the $L^2 (\Omega)$ boundedness to the case $|m -k| \leq \omega$ 
and the Schwarz inequality to the case $|m-k| > \omega$ 
for sequences $|m-k|^{-N}$, 
$|m-k|^{N} \| \chi_{C_\theta (m)} T \chi _{C_\theta (k)} \|_{L^2}$, 
and minimizing by taking suitable $\omega$. 
As to the second one \eqref{317-7}, we utilize the formula: 
$$
\psi ((M+\theta A)^{-1}) 
= (2\pi) ^{\frac{1}{2}} \int_{-\infty}^\infty e^{-it (M + \theta A)^{-1}} 
  \widehat {\psi } (t) \, dt .
$$
To estimate $\| e^{-it (M + \theta A)^{-1}}  \|_{A_N}$, 
we consider the commutator of 
$(x - \theta^{1/2} k)$ and $e^{-it (M + \theta A)^{-1}}$, 
which is rewritten with $\theta$, $(M+\theta A)^{-1}$, 
$\nabla(M + \theta A)^{-1}$ and is possible to be handled  
by the use of $L^2 (\Omega)$ boundedness, which proves \eqref{317-7}.  
As to the last one \eqref{317-8}, thanks to the formula 
$(M + \theta A)^{-\beta} 
= \frac{1}{\Gamma (\beta)} \int_0^\infty t^{\beta -1} e^{-Mt} e^{-t \theta A} dt$ 
and the Young inequality, we get
$$
\| (M + \theta A)^{-\beta} f \|_{\ell ^1 (L^2)_\theta} 
\leq \frac{1}{\Gamma (\beta)} \int_0^\infty t^{\beta -1} e^{-Mt} 
 \Big( \int_{\Omega} \| e^{-t\theta A} (\cdot ,y) \|_{\ell ^1 (L^2)_{\theta}}
    |f(y)| \, dy
 \Big) dt,
$$
where $\Gamma (\beta)$ is the Gamma function. 
By the Gaussian upper bounds of $e^{-t\theta A}$, we have 
$\| e^{-t\theta A} (\cdot ,y) \|_{\ell ^1 (L^2)_{\theta}} 
\leq C \theta^{-\frac{d}{4}} (1 + t ^{-\frac{d}{4}})$. 
These estimates yield the inequality \eqref{317-8}, 
since the integrability with respect to $t \in (0,\infty)$ is assured by 
$\beta > 4/d$.

\vskip3mm

\noindent 
{\bf Proof of Lemma \ref{lem:317-1}. } 
Instead of the inequality \eqref{317-1}, 
by replacing $2^{-j } \sqrt{A}, \sqrt{A}$ with $2^{-2j}A , A$, respectively, 
it is sufficient to show that 
\begin{equation}\label{317-2}
\| G(A)  \phi ( 2^{-2j} A) f \|_{L^p} 
\leq C \| G( 2^{2j} \cdot ) \phi (\cdot)
       \|_{H^{N+\frac{1}{2} + \delta} (\mathbb R)} \| f \|_{L^p}, 
\end{equation}
where ${\rm supp \, } \phi \subset [a^2 , b^2]$.  
Hence we show \eqref{317-2}. 

First we consider the case when $p = 1$. 
By decomposing $\Omega$ into cubes $C_{\theta} (k)$ and the H\"older inequality, 
we get 
\begin{equation}\label{317-3}
\| G(A) f \phi ( 2^{-2j} A) \|_{L^1} 
\leq C \theta ^{\frac{d}{2}} 
   \| G(A) \phi ( 2^{-2j} A) f \|_{\ell ^1 (L^2)_{\theta}} . 
\end{equation}
For fixed real numbers $M > 0$ and $ \beta > d/2$, let $\psi$ be such that 
\begin{equation}\label{317-10}
\psi (\mu) := G \big(2^{2j} (\mu^{-1} - M) \big) \phi (\mu ^{-1} -M) \mu ^{-\beta} . 
\end{equation}
It is easy to check that that 
$$
\psi \in C_0^\infty ((0,\infty)) 
\quad \text{and} \quad 
{\rm supp \, } \psi \subset 
\Big[ \frac{1}{M+b} , \frac{1}{M + a} \Big], 
$$
and 
$$
G(\lambda) \phi ( 2^{-2j} \lambda) 
= G(2^{2j} \cdot 2^{-2j}\lambda) \phi ( 2^{-2j} \lambda)  \mu^{-\beta} \cdot \mu ^\beta
= \psi (\mu) \mu ^\beta, 
$$
where $\lambda$ and $\mu$ are real numbers with 
$$
2^{-2j} \lambda = \mu ^{-1} - M. 
$$
The above equality yieds that 
\begin{equation}\label{317-4}
G( A) \phi ( 2^{-2j} A)  
= \psi \big( (M+2^{-2j }A) ^{-1} \big) (M + 2^{-2j} A)^{-\beta}. 
\end{equation}
Then it follows from \eqref{317-3}, \eqref{317-4} and the estimate \eqref{317-8} in Lemma~\ref{lem:317-2} 
that 
\begin{equation}\label{317-5}
\begin{split}
& 
\| G(A) \phi ( 2^{-2j} A) f \|_{L^1} 
\\
& 
\leq 
C \theta ^{\frac{d}{2}} 
   \big\| \psi \big( (M+2^{-2j }A) ^{-1} \big) (M + 2^{-2j} A)^{-\beta} f \big\|_{\ell ^1 (L^2)_{\theta}} . 
\\
& 
\leq 
C \theta ^{\frac{d}{2}} 
   \big\| \psi \big( (M+2^{-2j }A) ^{-1} \big) \big\|_{\ell^1 (L^1)_{\theta} \to \ell^1 (L^2)_{\theta}}
   \big\| (M + 2^{-2j} A)^{-\beta} \big\|_{L^1 \to \ell ^1 (L^2)_{\theta}} 
   \| f \|_{L^1}. 
\\
& 
\leq C   \big\| \psi \big( (M+2^{-2j }A) ^{-1} \big) \big\|_{\ell^1 (L^1)_{\theta} \to \ell^1 (L^2)_{\theta}}
   \| f \|_{L^1}. 
\end{split}
\end{equation}
By comparing the estimates  \eqref{317-2} and \eqref{317-5}, 
all we have to do is to show that 
\begin{equation}\label{317-9}
\| \psi \big( (M+2^{-2j }A) ^{-1} \big) \|_{\ell^1 (L^2)_{\theta} \to \ell^1 (L^2)_{\theta}} 
\leq C 
   \| G( 2^j \cdot ) \phi (\cdot)
       \|_{H^{N+\frac{1}{2} + \delta} (\mathbb R)} . 
\end{equation}
To apply the estimate \eqref{317-6}, we consider the operator norms 
$\| \cdot \|_{L^2 \to L^2}$ and $\| \cdot \|_{\mathscr A_N}$ 
of $\psi \big( (M+2^{-2j }A) ^{-1} \big)$. 
On the operator norm $\| \cdot \|_{L^2 \to L^2}$, we have from 
$N > d/2$ and 
the embedding $H^{N+\frac{1}{2} + \delta} (\mathbb R) \hookrightarrow L^\infty (\mathbb R)$ that 
\begin{equation}\notag
\| \psi \big( (M+2^{-2j }A) ^{-1} \big) \|_{L^2 \to L^2} 
\leq \| \psi \|_{L^\infty (\mathbb R)} 
\leq \| \psi \|_{H^{N+\frac{1}{2} + \delta} (\mathbb R)} 
\end{equation}
for any $\delta > 0$. 
As to $\| \psi \big( (M+2^{-2j }A) ^{-1} \big) \|_{\mathscr A_N}$, 
by applying the estimate \eqref{317-7} and the H\"older inequality, 
for any $\delta > 0$ there exists $C > 0$ such that 
\begin{equation}\notag 
\begin{split}
\| \psi ( (M + \theta A)^{-1} ) \|_{\mathscr A_N} 
& 
\leq C \theta ^{\frac{N}{2}} 
   \int_{-\infty}^\infty (1+ |\xi|^2)^{\frac{N}{2}} | \widehat \psi (\xi) | d\xi
\\
& 
\leq C \theta ^{\frac{N}{2}} 
   \big\| (1+ |\xi|^2)^{-\frac{1}{2} - \delta} \big\|_{L^2 (\mathbb R)}
   \big\| (1+ |\xi|^2)^{\frac{N}{2} + \frac{1}{2} + \delta} \widehat \psi \big\|_{L^2 (\mathbb R)}
\\
& 
\leq C \theta ^{\frac{N}{2}} 
     \| \psi \|_{H^{N+\frac{1}{2} + \delta } (\mathbb R)} . 
\end{split}
\end{equation}
Then we deduce from the above two estimates and \eqref{317-6} that 
\begin{equation}\notag 
\begin{split}
& 
\| \psi \big( (M+2^{-2j }A) ^{-1} \big) \|_{\ell^1 (L^2)_{\theta} \to \ell^1 (L^2)_{\theta}} 
\\
& \leq 
C \Big\{ \| \psi \|_{H^{N + \frac{1}{2} + \delta} (\mathbb R) } 
      + \theta ^{-\frac{d}{4}} 
        \big( \theta ^{\frac{N}{2}}\| \psi \|_{H^{N+\frac{1}{2} + \delta } (\mathbb R)} 
        \big)^{\frac{d}{2N}} 
        \big( \| \psi \|_{H^{N+\frac{1}{2} + \delta } (\mathbb R)} 
        \big) ^{1-\frac{d}{2N}}
  \Big\}
\\
& 
\leq C \| \psi \|_{H^{N+\frac{1}{2} + \delta} (\mathbb R)}  . 
\end{split}
\end{equation}
Since $\psi$ is defined by \eqref{317-10} and the support is bounded 
and away from the origin, we see from the change of variables by $\mu = (\lambda + M)^{-1}$ 
that 
$$
\| \psi (\cdot) \|_{H^{N+\frac{1}{2} + \delta } (\mathbb R)} 
\leq C \| G( 2^{2j} \cdot ) \phi (\cdot) 
       \|_{H^{N+\frac{1}{2} + \delta} (\mathbb R)}. 
$$
Hence the estimate \eqref{317-9} is obtained by the above two estimates, 
and the estimate \eqref{317-2} in the case when $p = 1$ is proved. 

We next consider the case when $p = \infty$. Since the dual space of $L^1 (\Omega)$ is $L^\infty (\Omega)$ 
and $C_0 ^\infty (\Omega)$ is dense in $L^1 (\Omega)$, the following holds: 
$$
\| G(A) \phi (2^{-j} A) f \|_{L^\infty} 
= \sup _{g \in C_{0}^{\infty} , \| g \|_{L^1} = 1} 
  \Big| \int _{\Omega} \big( G(A) \phi (2^{-j} A) f\big) \overline{g} \, dx 
  \Big| . 
$$
On the right member of the above equality, we have from the duality argument 
for the operator $G(A) \phi (2^{-j} A) $, the H\"older inequality and the estimate \eqref{317-2} 
with $p = 1$ that 
\begin{equation}\notag 
\begin{split}
\Big|  \int _{\Omega} \big( G(A) \phi (2^{-j} A) f\big) \, \overline{g} \, dx 
\Big|
& 
= \Big|  _{\mathcal X_0 '} \langle G(A) \phi (2^{-j} A) f, g \rangle_{\mathcal X_0}
  \Big|
\\
& 
= \Big| _{\mathcal X_0 '} \langle f, G(A) \phi (2^{-j} A) g \rangle_{\mathcal X_0}
  \Big|
\\
& 
= \Big| \int _{\Omega} f\, \overline{ G(A) \phi (2^{-j} A) g } \, dx 
  \Big|
\\
& 
\leq \| f \|_{L^\infty} \| G(A) \phi (2^{-j} A) g \|_{L^1}
\\
& 
\leq \| f \|_{L^\infty} \| G( 2^{2j} \cdot ) \phi (\cdot) 
       \|_{H^{N+\frac{1}{2} + \delta} (\mathbb R)} \| g \|_{L^1} ,
\end{split}
\end{equation}
where $g \in C_0 ^\infty$. 
These prove \eqref{317-2} in the case when $p = \infty$. 

As to the case when $1 < p < \infty$, 
the Riesz Thorin theorem allows us to obtain the estimate \eqref{317-2}. 
The proof of Lemma~\ref{lem:317-1} is complete. 
\hfill $\Box$

\section{Proof of Theorem \ref{thm:bdd}} 

\quad 
We prove Theorem \ref{thm:bdd} in this section. 
\\

\noindent 
{\bf Proof of \eqref{317-11}. } 
Put $\Phi _j := \phi_{j-1} + \phi_j + \phi_{j+1}$.
By applying the estimate \eqref{317-1} in Lemma~\ref{lem:317-1} with 
$$
G = G_t(\lambda) = e^{-t \lambda ^\alpha}, 
$$ 
we have 
\begin{equation}\notag 
\begin{split}
\big\| \phi_j (\sqrt{A}) e^{-t A^{\frac{\alpha}{2}}} f \big\|_{L^p} 
& 
= \big\| \big( G_t (\sqrt{A}) \Phi_j (\sqrt{A}) \big) \big( \phi_j (\sqrt{A} ) f\big)  \big\|_{L^p} 
\\
& 
\leq C \big\| G_t (2^j \sqrt{\cdot}\, ) \Phi_0 (\sqrt{\cdot} \, ) \big\|_{H^{N+\frac{1}{2}+\delta } (\mathbb R)} 
      \big\|  \phi_j (\sqrt{A}) f \big\| _{L^p} ,
\end{split}
\end{equation}
where $N > d/2$ and $\delta > 0$. 
Here it is easy to check that there exists $C > 0$ such that 
\begin{equation}\notag 
\big\| G_t (2^j \sqrt{\cdot}\,) \Phi_0 (\sqrt{\cdot}\,) \big\|_{H^{N+\frac{1}{2}+\delta} (\mathbb R)}  
\leq C e^{-C^{-1} t 2^{\alpha j} }
\quad \text{for any } j \in \mathbb Z, 
\end{equation}
and hence, 
\begin{equation}\label{317-13}
\begin{split}
\big\| \phi_j (\sqrt{A}) e^{-t A^{\frac{\alpha}{2}}} f \big\|_{L^p} 
\leq C e^{-C^{-1} t 2^{\alpha j} } 
     \big\| \phi_j (\sqrt{A}) f \big\|_{L^p}
     \quad \text{for any } j \in \mathbb Z .
\end{split}
\end{equation}
By multiplying $2^{sj}$ and taking the $\ell ^q (\mathbb Z)$ norm in the above inequality, 
we obtain the assertion \eqref{317-11}. 
\hfill $\Box$
\\

\noindent 
{\bf Proof of \eqref{317-12}. } 
By the inequalities 
$$
\big\| e^{-tA^{\frac{\alpha}{2}}} f \big\|_{\dot B^{s_2}_{p_2,q_2}} 
\leq \big\| e^{-tA^{\frac{\alpha}{2}}} f \big\|_{\dot B^{s_2}_{p_2,1}}, 
\quad 
\| f \|_{\dot B^{s_1}_{p_1,\infty}} 
\leq \| f \|_{\dot B^{s_1}_{p_1,q_1}}, 
$$
which are assured from the embedding relations in the Besov spaces, 
and taking $s_1 = 0 $ for the sake of simplicity, 
it is sufficient to show the following. 
\begin{equation}\label{317-14}
\big\| e^{-tA^{\frac{\alpha}{2}}} f \big\|_{\dot B^{s_2}_{p_2,1}}
\leq C t^{-\frac{d}{\alpha} (\frac{1}{p_1} - \frac{1}{p_2}) 
          - \frac{s_2}{\alpha} 
         } 
\| f \|_{\dot B^{0}_{p_1,\infty}} ,
\end{equation}
where 
$$
s_2 \geq 0,  
\quad p_1 \leq p_2 
\quad \text{and} \quad 
d \Big( \frac{1}{p_1} - \frac{1}{p_2} \Big) 
 + s_2 > 0 .
$$
We show \eqref{317-14}. 
It follows from the embedding $\dot B^{s_2 + d(\frac{1}{p_1} - \frac{1}{p_2})}_{p_1,1} 
\hookrightarrow \dot B^{s_2}_{p_2,1}$ and the estimate \eqref{317-13} that 
$$
\big\| e^{-tA^{\frac{\alpha}{2}}} f \big\|_{\dot B^{s_2}_{p_2,1}}
\leq C \big\| e^{-tA^{\frac{\alpha}{2}}} f \big\|_{\dot B^{s_2 + d(\frac{1}{p_1} - \frac{1}{p_2})}_{p_1,1}}
\leq C \sum _{ j \in \mathbb Z} 
   2^ {s_2 j + d(\frac{1}{p_1} - \frac{1}{p_2})j } e^{-ct 2^{\alpha j}}
   \| \phi_j (\sqrt{A}) f \|_{L^{p_1}} . 
$$
Since $s_ 2 + d (1/p_1 - 1/p_2) > 0$, we get 
\begin{equation}\notag 
\begin{split}
& 
\sum _{ j \in \mathbb Z} 
   2^ {s_2 j + d(\frac{1}{p_1} - \frac{1}{p_2})j } e^{-ct 2^{\alpha j}}
   \| \phi_j (\sqrt{A}) f \|_{L^{p_1}}
\\
& 
=
t^{-\frac{s_2}{\alpha} - \frac{d}{\alpha} (\frac{1}{p_1} - \frac{1}{p_2})}
 \sum _{ j \in \mathbb Z} 
\Big\{   (t 2^{\alpha j})^ { \frac{s_2}{\alpha} + \frac{d}{\alpha}(\frac{1}{p_1} - \frac{1}{p_2}) } e^{-ct 2^{\alpha j}}
\Big\} 
   \| \phi_j (\sqrt{A}) f \|_{L^{p_1}}
\\
& 
\leq C t^{-\frac{s_2}{\alpha} - \frac{d}{\alpha} (\frac{1}{p_1} - \frac{1}{p_2})}
 \| f \|_{\dot B^0_{p_1 ,\infty}}, 
\end{split}
\end{equation}
which proves \eqref{317-14}. Then the proof of \eqref{317-12} is complete. 
\hfill $\Box$

\section{Proof of Theorem \ref{thm:cont}} 

\quad 
We prove Theorem \ref{thm:cont} in this section. 
\\

\noindent 
{\bf Proof of {\rm (i)}. } Let $f \in \dot B^s_{p,q}(A)$. We take $f_N$ such that 
$$
f_N := \sum _{|j| \leq N} \phi_j (\sqrt{A}) f \quad \text{for } N \in \mathbb N .
$$
Since $q < \infty$, for any $\varepsilon > 0$ there exists $N_0 \in \mathbb N$ 
such that 
\begin{equation}\notag 
\| f_N - f \|_{\dot B^s_{p,q}} < \varepsilon 
\quad \text{for any } N \geq N _0
\end{equation}
The above inequality and boundedness \eqref{317-11}
 in Theorem~\ref{thm:bdd} imply that 
\begin{equation}\notag 
\begin{split}
\big\| e^{-t A^{\frac{\alpha}{2}}} f - f \big\|_{\dot B^s_{p,q}} 
& 
\leq  
\big\| e^{-t A^{\frac{\alpha}{2}}} f_N - f_N \big\|_{\dot B^s_{p,q}} 
+ 
\big\| e^{-t A^{\frac{\alpha}{2}}} (f_N -f) \big\|_{\dot B^s_{p,q}} 
+ 
\| f_N -f \|_{\dot B^s_{p,q}} 
\\
& 
\leq  
\big\| e^{-t A^{\frac{\alpha}{2}}} f_N - f_N \big\|_{\dot B^s_{p,q}} 
+ 
C \| f_N -f \|_{\dot B^s_{p,q}} 
\\
& 
\leq  
\big\| e^{-t A^{\frac{\alpha}{2}}} f_N - f_N \big\|_{\dot B^s_{p,q}} 
+ 
C \varepsilon 
\end{split}
\end{equation}
for any $t > 0$ provided that $N \geq N _0$. 
Then all we have to do is to show that 
\begin{equation}\label{317-15}
\lim _{ t \to 0} \big\| e^{-t A^{\frac{\alpha}{2}}} f_N - f_N \big\|_{\dot B^s_{p,q}} 
= 0. 
\end{equation}
We prove \eqref{317-15}. 
Noting that the spectrum of $f_N$ is restricted and 
$$
\big\| e^{-t A^{\frac{\alpha}{2}}} f_N - f_N \big\|_{\dot B^s_{p,q}} 
= \Big\{ \sum _{j = -N-1}^{N+1} \Big( 2^{sj}  
      \big\| \phi_j (\sqrt{A}) \big( e^{-t A^{\frac{\alpha}{2}}}-1 \big) f_N
      \big\|_{L^p} 
      \Big) ^q
  \Big\}^{\frac{1}{q}}, 
$$
we may consider the convergence of 
$\big\| \phi_j (\sqrt{A}) \big( e^{-t A^{\frac{\alpha}{2}}}-1 )f_N \big) 
 \big\|_{L^p} $ 
      for each $j$. 
For each $j = 0, \pm1, \pm2, \dots, \pm (N+1)$, 
it follows from \eqref{317-1} in Lemma \ref{lem:317-1} with 
$$
G =G_t(\lambda) = e^{-t \lambda ^{\alpha}}-1
$$
that 
\begin{equation}\notag 
\begin{split}
      \big\| \phi_j (\sqrt{A}) \big( e^{-t A^{\frac{\alpha}{2}}}-1 \big) f_N
      \big\|_{L^p} 
& 
=     \big\| \big( G_t (\sqrt{A} ) \Phi_j (\sqrt{A}) \big)
             \big( \phi_j (\sqrt{A}) f_N \big)
      \big\|_{L^p} 
\\
& 
\leq C 
     \big\| G_t ( 2^j \sqrt{\cdot}\, ) \Phi_0 ( \sqrt{\cdot}\, ) \big\|_{H^{N+\frac{d}{2}+\delta}}
      \big\| \phi_j (\sqrt{A}) f_N 
      \big\|_{L^p} , 
\end{split}
\end{equation}
where $\Phi_j := \phi_{j-1} + \phi_j + \phi_{j+1}$. 
Here it is readily checked that 
$$
\lim _{t \to 0 }  \big\| G_t ( 2^j \sqrt{\cdot}\, ) \Phi_0 ( \sqrt{\cdot}\, ) \big\|_{H^{N+\frac{d}{2}+\delta}}
= 0 
\quad \text{for each } j ,
$$
and hence, \eqref{317-15} is obtained. The proof of {\rm (i)} in Theorem \ref{thm:cont} is complete. 
\hfill $\Box$
\\

\noindent 
{\bf Proof of {\rm (ii)}. }
Put $\Phi_j := \phi_{j-1} + \phi_j + \phi_{j+1}$. 
By considering the dual operator of 
$e^{-tA^{ \frac{\alpha }{2} }} -1$, 
we have the following equality 
\begin{equation}\label{317-16}
\begin{split}
& 
 \sum _{ j \in \mathbb Z} 
\int _{\Omega} \Big\{ \phi_j (\sqrt{A}) \big( e^{-tA^{ \frac{\alpha }{2} }} -1 \big) f \Big\} 
     \overline{ \Phi_j (\sqrt{A}) g} \, dx  
\\
& 
= \sum _{ j \in \mathbb Z} 
\int _{\Omega} \Big\{ \phi_j (\sqrt{A})  f \Big\} 
     \overline{ \Phi_j (\sqrt{A}) \big( e^{-tA^{ \frac{\alpha }{2} }} -1 \big) g} \, dx .  
\end{split}
\end{equation}
It follows from the H\"older inequality that 
\begin{equation}
\begin{split}
& 
\sum _{ j \in \mathbb Z} 
\int _{\Omega} 
  \Big| \Big\{ \phi_j (\sqrt{A})  f \Big\} 
     \overline{ \Phi_j (\sqrt{A}) \big( e^{-tA^{ \frac{\alpha }{2} }} -1 \big) g} 
  \Big|\, dx
\\
& 
 \leq \sum _{j \in \mathbb Z} 2^{sj}
      \big\| \phi_j (\sqrt{A})  f \big\|_{L^p} 
      \cdot 2^{-sj}
      \big\| \Phi_j (\sqrt{A}) \big( e^{-tA^{ \frac{\alpha }{2} }} -1 \big) g \big\|_{L^{p'}}
\\
& 
\leq C \| f \|_{\dot B^s_{p,\infty}} 
     \big\| \big( e^{-tA^{ \frac{\alpha }{2} }} -1 \big) g \big\| _{\dot B^{-s}_{p',1}} , 
\end{split}
\end{equation}
which assures the absolute convergence of the series in \eqref{317-16} 
by the boundedness of $e^{-tA^{ \frac{\alpha }{2} }}$ in $\dot B^{-s}_{p',1}$ from \eqref{317-11} 
in Theorem~\ref{thm:bdd}. The above estimate and the assertion {\rm (i)} of Theorem~\ref{thm:cont} 
imply that 
\begin{equation} \notag 
\begin{split}
\Big| \sum _{ j \in \mathbb Z} 
\int _{\Omega} \Big\{ \phi_j (\sqrt{A}) \big( e^{-tA^{ \frac{\alpha }{2} }} -1 \big) f \Big\} 
     \overline{ \Phi_j (\sqrt{A}) g} \, dx  
\Big| 
& 
\leq C \| f \|_{\dot B^s_{p,\infty}} 
     \big\| \big( e^{-tA^{ \frac{\alpha }{2} }} -1 \big) g \big\| _{\dot B^{-s}_{p',1}}
\\
& 
\to 0 \quad \text{as  } t \to 0.
\end{split}
\end{equation}
This completes the proof of the assertion {\rm (ii)}. 
\hfill $\Box$

\section{Proof of Theorem \ref{thm:equiv}} 

\quad 
We prove Theorem \ref{thm:equiv} in this section. 
We first introduce the following lemma.

\begin{lem}\label{lem:317-3}
Let $\alpha > 0$, $s_0 \in \mathbb R$, $1 \leq p \leq \infty$. 
Then there exists $C > 0$ such that 
\begin{equation}\label{317-17}
\begin{split}
& 
C^{-1} (t 2^{\alpha j} ) ^{s_0} e^{-C t 2^{\alpha j} }
\big\| \phi_j (\sqrt{A}) f \big\|_{L^p}
\\
& 
\leq 
\big\| (tA ^{\frac{\alpha}{2}}) ^{s_0} e^{-t A^{\frac{\alpha}{2}}} \phi_j (\sqrt{A}) f 
\big\|_{L^p} 
\leq C ( t 2^{\alpha j} ) ^{s_0} e^{-C^{-1} t 2^{\alpha j} }
\big\| \phi_j (\sqrt{A}) f \big\|_{L^p}
\end{split}
\end{equation}
for any $t > 0$, $j \in \mathbb Z$ and $f \in L^p (\Omega)$. 
\end{lem}

\noindent 
{\bf Proof. } Put $\Phi _j := \phi_{j-1} + \phi_j + \phi_{j+1}$.
We start by proving the second inequality of the estimate  \eqref{317-17}. 
By applying the estimate \eqref{317-1} in Lemma~\ref{lem:317-1} with 
$$
G = G_t(\lambda) = (t \lambda ^\alpha)^{s_0} e^{-t \lambda ^\alpha}, 
$$ 
we have 
\begin{equation}\label{317-19}
\begin{split}
\big\| (tA ^{\frac{\alpha}{2}}) ^{s_0} e^{-t A^{\frac{\alpha}{2}}} \phi_j (\sqrt{A}) f 
\big\|_{L^p} 
& 
= 
 \big\| \big( G_t (\sqrt{A}) \Phi_j (\sqrt{A}) \big) \big( \phi_j (\sqrt{A} ) f\big)  \big\|_{L^p} 
\\
& 
\leq C \big\| G_t (2^j \sqrt{\cdot}\,) \Phi_0 (\sqrt{\cdot}\,) \big\|_{H^{N+\frac{1}{2}+\delta } (\mathbb R)} 
      \big\|  \phi_j (\sqrt{A}) f \big\| _{L^p} ,
\end{split}
\end{equation}
where $N > d/2$ and $\delta > 0$. 
Here it is easy to check that there exists $C > 0$ such that 
\begin{equation}\label{317-20}
\big\| G_t (2^j \sqrt{\cdot}\,) \Phi_0 (\sqrt{\cdot}\,) \big\|_{H^{N+\frac{1}{2}+\delta} (\mathbb R)}  
\leq C ( t 2^{\alpha j} ) ^{s_0} e^{-C^{-1} t 2^{\alpha j} } 
\quad \text{for any } j \in \mathbb Z, 
\end{equation}
and hence, 
\begin{equation}\notag 
\begin{split}
\big\| (tA ^{\frac{\alpha}{2}}) ^{s_0}  e^{-t A^{\frac{\alpha}{2}}} \phi_j (\sqrt{A}) f \big\|_{L^p} 
\leq C ( t 2^{\alpha j} ) ^{s_0} e^{-C^{-1} t 2^{\alpha j} } 
     \big\| \phi_j (\sqrt{A}) f \big\|_{L^p}
     \quad \text{for any } j \in \mathbb Z .
\end{split}
\end{equation}
This proves the second inequality of \eqref{317-17}. 

We turn to the first inequality of \eqref{317-17}. Since $\phi_j (\sqrt{A}) f$ is written as 
\begin{gather}\notag
\begin{split}
\phi_j (\sqrt{A}) f 
& 
= \Big( (t A^{\frac{\alpha}{2}}) ^{-s_0} e^{t A ^{\frac{\alpha}{2}}} \Phi_j (\sqrt{A}) 
  \Big)
  \Big( (t A^{\frac{\alpha}{2}}) ^{s_0 } e^{-t A^{\frac{\alpha}{2}}} \phi_j (\sqrt{A} ) f 
  \Big)
\\
& 
=: \Big( (t A^{\frac{\alpha}{2}}) ^{-s_0} e^{t A ^{\frac{\alpha}{2}}} \Phi_j (\sqrt{A}) 
  \Big)
  F, 
\end{split}
\end{gather}
all we have to do is to show that 
\begin{equation}\label{317-18} 
\big\| (t A^{\frac{\alpha}{2}}) ^{-s_0} e^{t A ^{\frac{\alpha}{2}}} \Phi_j (\sqrt{A})  F 
\big\|_{L^p} 
\leq C (t 2^{\alpha j}) ^{-s_0} e^{C t \lambda ^{\alpha}} \| F \|_{L^p} .
\end{equation}
Applying \eqref{317-1} in Lemma~\ref{lem:317-1} with 
$$
G = \tilde G_t (\lambda)  = (t \lambda ^\alpha)^{-s_0} e^{t \lambda ^\alpha}
$$
to the left member of \eqref{317-18}, we have from the similar argument to \eqref{317-19} and \eqref{317-20} 
that 
\begin{equation}\notag 
\begin{split}
\big\| (t A^{\frac{\alpha}{2}}) ^{-s_0} e^{t A ^{\frac{\alpha}{2}}} \Phi_j (\sqrt{A})  F 
\big\|_{L^p} 
& 
\leq C \| \tilde G_t (2^j \sqrt{\cdot}\,) \Phi_0 (\sqrt{\cdot}\,) \|_{H^{N + \frac{1}{2} + \varepsilon} (\mathbb R) } 
       \| F \|_{L^p}
\\
& 
\leq C (t 2^{\alpha j}) ^{-s_0} e^{C t \lambda ^{\alpha}} \| F \|_{L^p} .
\end{split}
\end{equation}
This proves \eqref{317-18} and the first inequality of \eqref{317-17} is obtained. 
Therefore we complete the proof of Lemma~\ref{lem:317-3}
\hfill $\Box$ 
\\

In what follows, we show the inequality \eqref{318-1} for $f \in \dot B^s_{p,q} (A)$ to 
prove Theorem \ref{thm:equiv}. 
 We note that the proof below concerns with the case when $q < \infty$ only, 
since the case when $q = \infty$ is also shown analogously with some 
modification.
\\

\noindent 
{\bf Proof of the first inequality of \eqref{318-1}. } 
By the embedding $L^p (\Omega), \dot B^0_{p, r} (A) \hookrightarrow \dot B^0_{p,\infty} (A)$, 
it is sufficient to show that 
\begin{equation}\label{317-21}
C^{-1}\| f \|_{\dot B^s_{p,q}} 
\leq 
\Big\{ \int _0^\infty 
     \Big(  t^{-\frac{s}{\alpha}} 
            \| (tA^{\frac{\alpha}{2}} ) ^{s_0} e^{-tA^{\frac{\alpha}{2}}} f \|_{\dot B^0_{p,\infty}} 
     \Big)^q
     \frac{dt}{t}
\Big\} ^{\frac{1}{q}}. 
\end{equation}
We have from the definition of norm $\| \cdot \|_{\dot B^0_{p,\infty} }$ 
and the first inequality of estimate \eqref{317-17} in Lemma \ref{lem:317-3} 
that 
\begin{equation}\notag 
\begin{split}
& 
\Big\{ \int _0^\infty 
     \Big(  t^{-\frac{s}{\alpha}} 
            \| (tA^{\frac{\alpha}{2}} ) ^{s_0} e^{-tA^{\frac{\alpha}{2}}} f \|_{\dot B^0_{p,\infty}} 
     \Big)^q
     \frac{dt}{t}
\Big\} ^{\frac{1}{q}}
\\
& 
\geq C^{-1} 
\Big\{ \int _0^\infty 
     \Big(  t^{-\frac{s}{\alpha}} 
           \sup _{ j \in \mathbb Z} 
             (t 2^{\alpha j}) ^{s_0} e^{-C t2^{\alpha j}}
             \|  \phi_j (\sqrt{A}) f \|_{L^p} 
     \Big)^q
     \frac{dt}{t}
\Big\} ^{\frac{1}{q}} . 
\end{split}
\end{equation}
Decomposing $(0,\infty)$ in the last line by 
\begin{equation}\label{317-24}
(0,\infty) = \bigcup_{k \in \mathbb Z} [ 2^{-\alpha (k+1)} , 2^{-\alpha k}], 
\end{equation}
we get 
\begin{equation}\label{317-22}
\begin{split}
& 
\Big\{ \int _0^\infty 
     \Big(  t^{-\frac{s}{\alpha}} 
            \| (tA^{\frac{\alpha}{2}} ) ^{s_0} e^{-tA^{\frac{\alpha}{2}}} f \|_{\dot B^0_{p,\infty}} 
     \Big)^q
     \frac{dt}{t}
\Big\} ^{\frac{1}{q}}
\\
& 
\geq C^{-1} 
\Big\{ \sum _{ k \in \mathbb Z} \int _{2^{-\alpha (k+1)}} ^{2^{-\alpha k}}
     \Big(  t^{-\frac{s}{\alpha}} 
           \sup _{ j \in \mathbb Z} 
             (t 2^{\alpha j}) ^{s_0} e^{-C t2^{\alpha j}}
             \|  \phi_j (\sqrt{A}) f \|_{L^p} 
     \Big)^q
     \frac{dt}{t}
\Big\} ^{\frac{1}{q}}
\\
& 
\geq C^{-1} 
\Big\{ \sum _{ k \in \mathbb Z}
     \Big(  2^{sk} 
           \sup _{ j \in \mathbb Z} 
             ( 2^{\alpha (j-k)}) ^{s_0} e^{-C 2^{\alpha (j-k)}}
             \|  \phi_j (\sqrt{A}) f \|_{L^p} 
     \Big)^q
\Big\} ^{\frac{1}{q}} . 
\end{split}
\end{equation}
Here it follows from the H\"older inequality that 
\begin{equation}\notag 
\begin{split}
&           \sup _{ j \in \mathbb Z} 
             ( 2^{\alpha (j-k)}) ^{s_0} e^{-C 2^{\alpha (j-k)}}
             \|  \phi_j (\sqrt{A}) f \|_{L^p} 
\\
& 
\geq C^{-1} \Big\{ 
      \sum _{j \in \mathbb Z} 
      \Big( \frac{1}{1+ \alpha ^2 |j-k|^2} 
           \cdot ( 2^{\alpha (j-k)}) ^{s_0} e^{-C 2^{\alpha (j-k)}}
             \|  \phi_j (\sqrt{A}) f \|_{L^p} 
      \Big)^q
      \Big\}^{\frac{1}{q}} . 
\end{split}
\end{equation}
Then we deduce from \eqref{317-22} and the above inequality that 
\begin{equation}\notag 
\begin{split}
& 
\Big\{ \int _0^\infty 
     \Big(  t^{-\frac{s}{\alpha}} 
            \| (tA^{\frac{\alpha}{2}} ) ^{s_0} e^{-tA^{\frac{\alpha}{2}}} f \|_{\dot B^0_{p,\infty}} 
     \Big)^q
     \frac{dt}{t}
\Big\} ^{\frac{1}{q}}
\\
& 
\geq C^{-1} 
\Big\{ \sum _{ k \in \mathbb Z}
      \big( 2^{sk} \big) ^q
           \sum _{ j \in \mathbb Z} 
      \Big( \frac{1}{1+ \alpha ^2 |j-k|^2} 
           \cdot ( 2^{\alpha (j-k)}) ^{s_0} e^{-C 2^{\alpha (j-k)}}
             \|  \phi_j (\sqrt{A}) f \|_{L^p} 
      \Big)^q
\Big\} ^{\frac{1}{q}}
\\
& 
= C^{-1} 
\Big\{ \sum _{ j \in \mathbb Z} 
         \big( 2^{sj} \|  \phi_j (\sqrt{A}) f \|_{L^p} \big)^q
       \sum _{ k \in \mathbb Z}
      \Big( \frac{2^{-s(j-k)}}{1+ \alpha ^2 |j-k|^2} 
           \cdot ( 2^{\alpha (j-k)}) ^{s_0} e^{-C 2^{\alpha (j-k)}}
      \Big)^q
\Big\} ^{\frac{1}{q}}
\\
& 
= C^{-1} \| f \|_{\dot B^s_{p,q}}
\Big\{ \sum _{ k \in \mathbb Z}
      \Big( \frac{2^{(s_0 \alpha - s)k}}{1+ \alpha ^2 |k|^2} 
           \cdot  e^{-C 2^{\alpha k}}
      \Big)^q
\Big\} ^{\frac{1}{q}}.
\end{split}
\end{equation}
Since $s_0 > s/\alpha$ and the summation appearing in the last line 
converges, we obtain \eqref{317-21}. 
Hence the proof of the first inequality of \eqref{318-1} is complete. 
\hfill $\Box$
\\

\noindent 
{\bf Proof of the second inequality of \eqref{318-1}. }
By the embedding $\dot B^0_{p,1} (A) \hookrightarrow L^p (\Omega), \dot B^0_{p, q} (A)$, 
it is sufficient to show that 
\begin{equation}\label{317-23}
\Big\{ \int _0^\infty 
     \Big(  t^{-\frac{s}{\alpha}} 
            \| (tA^{\frac{\alpha}{2}} ) ^{s_0} e^{-tA^{\frac{\alpha}{2}}} f \|_{\dot B^0_{p,1}} 
     \Big)^q
     \frac{dt}{t}
\Big\} ^{\frac{1}{q}} 
\leq C \| f \|_{\dot B^s_{p,q} (A)}. 
\end{equation}
Analogously to the proof of \eqref{317-21}, 
we apply the second inequality of \eqref{317-17} in Lemma \ref{lem:317-3} 
instead of the first one and the decomposition \eqref{317-24} to get 
\begin{equation}\notag 
\begin{split}
& 
\Big\{ \int _0^\infty 
     \Big(  t^{-\frac{s}{\alpha}} 
            \| (tA^{\frac{\alpha}{2}} ) ^{s_0} e^{-tA^{\frac{\alpha}{2}}} f \|_{\dot B^0_{p,\infty}} 
     \Big)^q
     \frac{dt}{t}
\Big\} ^{\frac{1}{q}}
\\
& 
\leq \Big\{ \sum _{ k \in \mathbb Z}
     \Big(  2^{sk} 
           \sum _{ j \in \mathbb Z} 
             ( 2^{\alpha (j-k)}) ^{s_0} e^{-C^{-1} 2^{\alpha (j-k)}}
             \|  \phi_j (\sqrt{A}) f \|_{L^p} 
     \Big)^q
\Big\} ^{\frac{1}{q}} . 
\end{split}
\end{equation}
Here the H\"older inequality yields that 
\begin{equation}\notag 
\begin{split}
& \sum _{ j \in \mathbb Z} 
             ( 2^{\alpha (j-k)}) ^{s_0} e^{-C^{-1} 2^{\alpha (j-k)}}
             \|  \phi_j (\sqrt{A}) f \|_{L^p} 
\\
& \leq C \Big\{ \sum _{ j \in \mathbb Z} 
           \Big( (1+ \alpha^2 |j-k|^2)
             ( 2^{\alpha (j-k)}) ^{s_0} e^{-C^{-1} 2^{\alpha (j-k)}}
             \|  \phi_j (\sqrt{A}) f \|_{L^p} 
           \Big) ^q
       \Big\} ^{\frac{1}{q}}. 
\end{split}
\end{equation}
Then we have from the above two estimates that 
\begin{equation}\notag 
\begin{split}
& 
\Big\{ \int _0^\infty 
     \Big(  t^{-\frac{s}{\alpha}} 
            \| (tA^{\frac{\alpha}{2}} ) ^{s_0} e^{-tA^{\frac{\alpha}{2}}} f \|_{\dot B^0_{p,\infty}} 
     \Big)^q
     \frac{dt}{t}
\Big\} ^{\frac{1}{q}}
\\
& 
\leq C \Big\{ \sum _{ k \in \mathbb Z}
     \big(  2^{sk} \big) ^q
           \sum _{ j \in \mathbb Z} 
           \Big( (1+ \alpha^2 |j-k|^2)
             ( 2^{\alpha (j-k)}) ^{s_0} e^{-C^{-1} 2^{\alpha (j-k)}}
             \|  \phi_j (\sqrt{A}) f \|_{L^p} 
           \Big) ^q
\Big\} ^{\frac{1}{q}}
\\
& 
= C \Big\{  \sum _{ j \in \mathbb Z} 
           \big( 2^{sj} \|  \phi_j (\sqrt{A}) f \|_{L^p}  \big) ^q 
            \sum _{ k \in \mathbb Z}
           \Big( 2^{-s (j-k)}(1+ \alpha^2 |j-k|^2)
             ( 2^{\alpha (j-k)}) ^{s_0} e^{-C^{-1} 2^{\alpha (j-k)}}
           \Big) ^q
\Big\} ^{\frac{1}{q}}
\\
& 
= C \| f \|_{\dot B^s_{p,q}} 
\Big\{      \sum _{ k \in \mathbb Z}
           \Big( (1+ \alpha^2 |k|^2)
             2^{ (s_0 \alpha -s) k}  e^{-C^{-1} 2^{\alpha k}}
           \Big) ^q
\Big\} ^{\frac{1}{q}} .
\end{split}
\end{equation}
Since $s_0 > s/\alpha$ and the summation appearing in the last line 
converges, we obtain \eqref{317-23}. 
The proof of the second inequality of \eqref{318-1} is complete. 
\hfill $\Box$

\section{Proof of Theorem \ref{thm:maxi}}

\quad 
We prove Theorem \ref{thm:maxi} in this section. 
\\

\noindent 
{\bf Proof of \eqref{318-2}. } 
It is sufficient to prove the case when $s = 0$ 
thanks to the 
lifting property in the proposition in section 1. 
We also consider the case when $q < \infty$ only, 
since the case when $q = \infty$ is also shown analogously.
First we prove that 
\begin{equation}\label{317-24-2}
\| A^{\frac{\alpha}{2}} u \|_{L^q (0,\infty ; \dot B^0_{p,q})} 
\leq 
C\| u_0 \|_{\dot B^{\alpha - \frac{\alpha}{q}}_{p,q}  }
+ C \| f \| _{L^q (0,\infty ; \dot B^0_{p,q})} . 
\end{equation}
By the definition of $u$ and the triangle inequality, we get 
\begin{equation}\label{317-25}
\begin{split}
& 
\| A^{\frac{\alpha}{2}} u \|_{L^q (0,\infty ; \dot B^0_{p,q})} 
\\
& 
\leq \| A^{\frac{\alpha}{2}} e^{- t A^{\frac{\alpha}{2}}} u_0 \|_{L^q (0,\infty ; \dot B^0_{p,q})} 
  + \Big\| A^{\frac{\alpha}{2}} \int_0^t e^{ -(t-\tau) A^{\frac{\alpha}{2}}} f(\tau) \, d\tau 
    \Big\|_{L^q (0,\infty ; \dot B^0_{p,q})} . 
\end{split}
\end{equation}
On the first term of the right member in the above inequality, 
it follows from the estimate \eqref{318-1} for $s_0 = 1$, $s = \alpha -\alpha /q$ that 
\begin{equation}\label{317-27}
\| A^{\frac{\alpha}{2}} e^{- t A^{\frac{\alpha}{2}}} u_0 \|_{L^q (0,\infty ; \dot B^0_{p,q})} 
\leq C \| u_0 \| _{L^q (0,\infty ; \dot B^{\alpha - \frac{\alpha}{q}}_{p,q})}.
\end{equation}
As to the second one, we start by proving that 
\begin{equation}\label{317-26}
\begin{split}
& 
\Big\| \phi_j (\sqrt{A}) A^{\frac{\alpha}{2}} \int_0^t e^{ -(t-\tau) A^{\frac{\alpha}{2}}} f(\tau) \, d\tau 
\Big\|_{L^p} 
\\
& 
\leq C 2^{\frac{\alpha}{q} j} 
   \Big\{ \int_0^t 
   \Big( e^{-C^{-1} (t-\tau) 2^{\alpha j}} \| \phi_j (\sqrt{A}) f \|_{L^p} 
   \Big) ^q 
   \, d\tau
   \Big\}^{\frac{1}{q}} .
\end{split}
\end{equation}
The above estimate \eqref{317-26} is verified by applying 
the estimate \eqref{317-17} in Lemma \ref{lem:317-3} and the H\"older inequality, in fact, we get 
\begin{equation}\notag 
\begin{split}
& 
\Big\| \phi_j (\sqrt{A}) A^{\frac{\alpha}{2}} \int_0^t e^{ -(t-\tau) A^{\frac{\alpha}{2}}} f(\tau) \, d\tau 
\Big\|_{L^p} 
\\
& 
\leq C 
2^{\alpha j} 
  \int_0^t e^{ -C^{-1}(t-\tau) 2^{\alpha j}} 
       \| \phi_j (\sqrt{A}) f(\tau) \|_{L^p} \, d\tau 
\\
& 
\leq C 
2^{\alpha j} 
\| e^{- (2C)^{-1} (t-\tau) 2^{\alpha j}} \|_{L^{\frac{q}{q-1}} (\{ 0 \leq \tau \leq t\})}
  \Big\{ \int_0^t 
     \Big( e^{ -(2C)^{-1}(t-\tau) 2^{\alpha j}} 
       \| \phi_j (\sqrt{A}) f(\tau) \|_{L^p}
     \Big) ^q \, d\tau 
  \Big\} ^{\frac{1}{q}}
\\
& 
\leq C 
2^{\frac{\alpha}{q} j} 
  \Big\{ \int_0^t 
     \Big( e^{ -(2C)^{-1}(t-\tau) 2^{\alpha j}} 
       \| \phi_j (\sqrt{A}) f(\tau) \|_{L^p}
     \Big) ^q \, d\tau 
  \Big\} ^{\frac{1}{q}}. 
\end{split}
\end{equation}
By the estimate \eqref{317-26}, we have 
\begin{equation}\label{317-28}
\begin{split}
& \Big\| A^{\frac{\alpha}{2}} \int_0^t e^{ -(t-\tau) A^{\frac{\alpha}{2}}} f(\tau) \, d\tau 
    \Big\|_{L^q (0,\infty ; \dot B^0_{p,q})} 
\\
& \leq C 
 \bigg[ 
   \int_0 ^\infty\sum _{ j \in \mathbb Z} 
 \bigg\{ 2^{\frac{\alpha}{q} j} 
  \bigg( \int_0^t 
     \big( e^{ -(2C)^{-1}(t-\tau) 2^{\alpha j}} 
       \| \phi_j (\sqrt{A}) f(\tau) \|_{L^p}
     \big) ^q \, d\tau 
  \bigg) ^{\frac{1}{q}}
  \bigg\} ^q
  \,dt 
  \bigg] ^{\frac{1}{q}}
\\
& = C 
 \bigg[ 
   \int_0 ^\infty\sum _{ j \in \mathbb Z}  \| \phi_j (\sqrt{A}) f(\tau) \|_{L^p} ^q
      \bigg( 2^{\alpha j}  \int_\tau ^\infty e^{ -q(2C)^{-1}(t-\tau) 2^{\alpha j}} \,dt 
      \bigg)
      \, d\tau
  \bigg] ^{\frac{1}{q}}
\\
& 
= C \| f \| _{L^q (0,\infty ; \dot B^0_{p,q}(A))} . 
\end{split}
\end{equation}
Then the estimates \eqref{317-25}, \eqref{317-27} and \eqref{317-28} imply 
the inequality \eqref{317-24-2}. 
The estimate for $\partial _t u$, i.e., the inequality  
$$
\| \partial _t u \|_{L^q (0,\infty ; \dot B^0_{p,q})} 
\leq 
C\| u_0 \|_{\dot B^{\alpha - \frac{\alpha}{q}}_{p,q}  }
+ C \| f \| _{L^q (0,\infty ; \dot B^0_{p,q})} 
$$
is verified by the estimate \eqref{317-24-2} and the equality 
$$
\partial _t u = - A^{\frac{\alpha}{2} } u + f . 
$$
Hence we obtain the estimate \eqref{318-2} 
and the proof is complete. 
\hfill $\Box$

\section{Results for the inhomogeneous Besov spaces}

\quad 
We should mention that similar theorems 
also hold for the heat semi-group in the inhomogeneous Besov spaces 
$B^s_{p,q} (A)$. 
We also note that semi-group generated by the fractional Laplacian 
can not be treated analogously by the direct application of boundedness 
of scaled spectral multiplier in Lemma~\ref{lem:317-1} 
(see the comment below Theorem~\ref{thm:inhom})
\\

First we recall the definition of 
$B^s_{p,q} (A)$. Let $\psi$ is in $C_0 ^\infty ((-\infty,\infty)) $ such that 
$$
\psi (\lambda ^2) + \sum _{ j \in \mathbb N} \phi _j (\lambda) = 1 
\quad \text{for any } \lambda \geq 0 . 
$$
The inhomogeneous Besov spaces $B^s_{p,q} (A)$ is defined as 
follows (see \cite{IMT-preprint2}). 
\\

\noindent 
{\bf Definition. }
{\it 
For $s \in \mathbb R$ and $1 \leq p,q \leq \infty$, 
 $B^s_{p,q} (A) $ is defined by letting 
\begin{equation}\notag
B^s_{p,q} (A)
:= \{ f \in \mathcal X'_0 (\Omega) 
     \, | \, 
     \| f \|_{B^s_{p,q} (A)} < \infty
   \} , 
\end{equation}
where
\begin{equation}\notag 
     \| f \|_{B^s_{p,q} (A)} 
       := \| \psi ( A) f \|_{L^p  } 
          + \big\| \big\{ 2^{sj} \| \phi_j (\sqrt{A}) f \|_{L^p  }  
                 \big\}_{j \in \mathbb  N}
          \big\|_{\ell^q (\mathbb N)}. 
\end{equation}
}

\vskip3mm 

The high frequency part is able to be treated in the same way as the proof for 
the homogeneous case by using Lemma~\ref{lem:317-1}. 
As to the low frequency part, we employ the pointwise estimate of the kernel of 
$e^{-tA}$ 
$$
0 \leq e^{-tA} (x,y) 
\leq (4\pi t)^{-\frac{d}{2}} 
    \exp \Big( \frac{|x-y|^2}{4t} \Big),
$$
which assures the boundedness of $e^{-t A}$ in $L^p (\Omega)$ 
and also $B^s_{p,q} (A)$ as well as the case when $\Omega = \mathbb R^d$. 
In order to treat continuity in time of $e^{-tA}$, we need the following obtained 
by the similar proof to that of Lemma~\ref{lem:317-1}.

\begin{lem}\label{lem:530-1}
Let $N > d/2$, $1 \leq p \leq \infty$, $\delta > 0$, 
$\psi \in C_0 ^\infty ((-\infty , \infty))$ and 
$G \in H^{N + \frac{1}{2} + \delta} (\mathbb R)$. 
Then there exists a positive constant $C $ such that 
for any $f \in L^p (\Omega)$ 
\begin{gather}\label{405-1}
\| G(A) \psi ( A)  f \|_{L^p} 
\leq C \| G( \cdot ) \psi (\cdot)
       \|_{H^{N+\frac{1}{2} + \delta} (\mathbb R)} \| f \|_{L^p}
\end{gather}
\end{lem}

We take $G$ such that 
$$
G(\lambda) := 
e^{-t \lambda} -1 
\quad \text{for any } \lambda \in \mathbb R 
$$
to apply the above lemma. 
For the above $G$ it is easy to check that 
$$
\| G( \cdot ) \psi (\cdot)
       \|_{H^{N+\frac{1}{2} + \delta} (\mathbb R)} 
    \to 0 \quad \text{as } t \to 0.
$$
Hence for any $f \in B^s_{p,q} (A)$, it follows from \eqref{405-1} that 
$$
\lim _{ t\to 0} \| \psi (A) (e^{-tA} f -f) \|_{L^p} = 0 . 
$$
According to the boundedness and the continuity of $e^{-t A}$, 
we obtain the following result for the inhomogeneous Besov spaces. 

\begin{thm}\label{thm:inhom}
Let $s \in \mathbb R$, $1 \leq p,p_1,p_2,q \leq \infty$ and 
$1/p + 1/p' = 1$. 
Let $\Psi$ and $\Psi _j$ with $j \in \mathbb N$  
be such that 
\begin{gather}\notag 
\Psi(A) := \psi (A) + \phi_1 (\sqrt{A}), \quad 
\Phi_1 (\sqrt{A}):= \psi(A) + \phi_1 (\sqrt{A}) + \phi_2 (\sqrt{A}) , \quad 
\\
\notag 
\Phi_j (\sqrt{A}):= \phi_{j-1} (\sqrt{A}) + \phi_j(\sqrt{A}) + \phi_{j+1} (\sqrt{A})
\quad \text{for } j \geq 2.
\end{gather}
{\rm (i)} 
There exists a constant $C >0$ such that 
\begin{equation}\notag 
\big\| e^{-tA } f \big\|_{ B^{s}_{p,q}(A)} 
\leq C \| f \|_{ B^{s}_{p , q} (A) } 
\end{equation}
for any $f \in B^s_{p,q} (A)$. If $p_1 \leq p_2$, then there exists a constant $C > 0$ such that 
\begin{equation}\notag 
\big\| e^{-tA} f \big\|_{B^{s}_{p_2,q}(A)} 
\leq C t^{-\frac{d}{2} (\frac{1}{p_1} - \frac{1}{p_2}) }
     \| f \|_{B^{s}_{p_1 , q} (A) } 
\end{equation}
for any $f \in B^s_{p_1 ,q} (A)$. 
\\
{\rm (ii)} If $q < \infty$ and $f \in B^s_{p,q} (A)$, then 
\begin{equation}\notag 
\lim _{t \to 0} \big\| e^{-tA } f -f 
\big\|_{B^{s}_{p,q}(A)} = 0.
\end{equation}
If $q = \infty$, $1 < p \leq \infty$ and $f \in B^s_{p,\infty} (A)$, then 
$e^{-tA} f$ converges to $f$ in the dual weak sense as $t \to 0$, 
namely,  
\begin{equation}\notag 
\begin{split}
\lim _{ t \to 0 }
\bigg[ \int _{\Omega}
& 
  \Big\{ \psi (A)\big( e^{-tA } f -f \big) 
  \Big\} \, 
  \overline{ \Psi (A) g} \, dx 
\\
& 
+ 
\sum _{ j \in \mathbb N} 
\int _{\Omega} 
   \Big\{ \phi_j (\sqrt{A})\big( e^{-tA } f -f \big) 
   \Big\} \, 
   \overline{ \Phi_j (\sqrt{A}) g} \, dx 
\bigg]= 0
\end{split}
\end{equation}
for any $g \in \dot B^{-s}_{p',1} (A)$. 
\\
{\rm (iii)} Let $T>0$, $s, s_0 \in \mathbb R$, $s_0 > s/ 2$. Then 
\begin{equation}\notag
\| f \|_{B^s_{p,q} (A)} 
\simeq 
\| \psi (T A) f \|_{L^p} +
\Big\{ \int _0^T
     \Big(  t^{-\frac{s}{2}} \| (tA ) ^{s_0} e^{-tA} f \|_{X} 
     \Big)^q
     \frac{dt}{t}
\Big\} ^{\frac{1}{q}}
\end{equation}
for any $f \in B^s_{p,q} (A)$, 
where $X = L^p (\Omega), B^0_{p,r} (A)$ with $1 \leq r \leq \infty$. 
\\
{\rm (iv)} 
Let $T > 0$, $u_0 \in B^{s+ 2 - \frac{2}{q}}_{p,q} (A) $ and  
$ f \in L^q (0,T ; B^s_{p,q}(A))$. Assume that $u $ satisfies 
\begin{gather}
\notag 
u(t) = e^{- t A } u_0 
   + \int_0^t e^{-(t-\tau) A} f(\tau) d\tau . 
\end{gather}
Then there exists a constant $C = C(T)>0$ indeoendent of $u_0$ and $f$ such that 
\begin{equation}\notag 
\| \partial _t u \|_{L^q (0,T ; B^s_{p,q}(A))}
+ \| A^{\frac{\alpha}{2}} u \|_{L^q (0,T ; B^s_{p,q}(A))} 
\leq 
C\| u_0 \|_{B^{s+ 2 - \frac{2}{q}}_{p,q} (A) }
+ C \| f \| _{L^q (0,T ; B^s_{p,q}(A))} . 
\end{equation}

\end{thm}

\noindent 
{\bf Remark}. 
Let us mention what is obtained by the abstract theory 
for sectorial operators by Da Prato-Grisvard~\cite{DaGr-1975} 
(see also \cite{Hass_2006,Lu_1995}). 
Let $X = B^0_{p,q} (A)$. We can consider $A$ as 
a sectorial operator with the domain $D(A^{\alpha}) = B^{2}_{p,q} (A)$. 
Let $0 < T < \infty$, $1 < q < \infty$, $1 \leq p,r \leq \infty$, 
$\theta \in (0,1)$ and $\alpha > 0$. 
Then for any $f \in L^q (0,T ; (X , D(A^\alpha) )_{\theta,r})$
the equation 
\begin{equation}\notag 
\begin{cases}
\displaystyle \frac{du}{dt} + Au = f , 
& \quad 0 < t < T ,
\\
u(0) = 0
\end{cases}
\end{equation}
admits a unique solution $u$ satisfying 
$$
\Big\| \frac{du}{dt} \Big\|_{L^q (0,T; (X, D(A^\alpha))_{\theta, r} )}
+ \| A u \|_{L^q (0,T ; (X , D(A^\alpha) )_{\theta,r})} 
\leq C \| f \|_{L^q (0,T ; (X , D(A^\alpha) )_{\theta,r})} , 
$$
where $C$ depends on $T$. 
Here we note that 
$(X , D(A^\alpha) )_{\theta,r} = B^{2\alpha \theta}_{p,r} (A)$ 
and $2 \alpha \theta$ is possible to be an arbitrary positive number 
since $\alpha > 0$ and $\theta \in (0,1)$.

\vskip3mm

Let us give a few remark on the semi-group generated by $A^{\frac{\alpha}{2}}$. 
If we consider to apply Lemma \ref{lem:530-1} directly, 
it is impossible to obtain the boundedness of $e^{-t A^{\frac{\alpha}{2}}}$ 
for general $\alpha$. In fact, taking 
$$
G = G_t(\lambda) = e^{- t |\lambda|^{\frac{\alpha}{2}}}, 
$$
and applying \eqref{405-1}, we see that the $H^{N+\frac{1}{2}+\delta} (\mathbb R)$ norm 
of the above $G = G_t (\lambda)$ is not finite for small $\lambda > 0$ 
 because of less regularity around $\lambda = 0$. 
On the other hand, if $\alpha$ is even or sufficiently large, 
the $H^{N+\frac{1}{2}+\delta} (\mathbb R)$ norm of $e^{- t |\lambda|^{\frac{\alpha}{2}}}$ 
is finite and we can get some results. However this argument does not 
reach at the optimal estimate, and hence, we do not treat in this paper 
and will treat in the future work.

\appendix
\section{Real interpolation}

\quad 
In this appendix, we give a remark 
that real interpolation can be considered 
in the Besov spaces $\dot B^s_{p,q} (A)$ and $B^s_{p,q} (A)$ 
on open sets as well as the whole space case. 
We recall the definition of real interpolation spaces $(X_0, X_1)_{\theta,q}$ 
for Banach spaces $X_0$ and $X_1$ (see e.g. \cite{BL_1976,Pee_1968,Triebel_1983}). 
\\

\noindent 
{\bf Definition. }  
{\it 
Let $0 < \theta < 1$ and $ 1 \leq q \leq \infty $. 
$(X_0, X_1)_{\theta,q}$ is defined by letting 
$$
(X_0,X_1)_{\theta,q} 
:= \Big\{ a \in X_0 + X_1 \, \Big| \, 
      \| a \|_{(X_0,X_1)_{\theta,q}} := 
      \Big\{ \int_0^\infty \big( t^{-\theta} K(t,a) \big) ^q 
             \, \frac{dt}{t}
      \Big\}^{\frac{1}{q}}
      < \infty
   \Big\},
$$
where $K(t,a)$ is Peetre's K-function 
$$
K(t,a) := 
 \inf 
 \big\{ \| a_0 \|_{X_0} + t \| a_1 \|_{X_1} \,\big| \, 
       a = a_0 + a_1, \, a_0 \in X_0 , \, a_1 \in X_1 
 \big\}. 
$$
}

As well as the case when $\Omega = \mathbb R^d$, 
we obtain the following. 

\begin{prop}\label{prop:real}
Let $0 < \theta < 1$, $s,s_0,s_1 \in \mathbb R$ and $1 \leq p,q,q_0,q_1 \leq \infty$. 
Assume that $s_0 \not = s_1$ and 
$s = (1-\theta)s_0 + \theta s_1$. Then 
\begin{gather}
\notag 
\big(  \dot B^{s_0}_{p,q_0} (A)  , \, \dot B^{s_1}_{p,q_1} (A) \big)_{\theta , q} 
= \dot B^s_{p,q} (A), 
\\ \notag 
\big(  B^{s_0}_{p,q_0} (A)  , \,  B^{s_1}_{p,q_1} (A) \big)_{\theta , q} 
= B^s_{p,q} (A). 
\end{gather}
\end{prop}

We omit the proof of the above proposition since one can show 
analogously to the whole space case 
(see e.g. \cite{Triebel_1983}).

\vskip10mm 

\noindent
{\bf Acknowledgements.}
The author would like to thank the referee for his important comments. 
The author was supported by the Grant-in-Aid for Young Scientists (B) (No.~25800069)
from JSPS 
and by JSPS Program for Advancing Strategic International Networks to Accelerate 
the Circulation of Talented Researchers.

\end{document}